\documentclass[3p]{elsarticle}

\usepackage{amsmath}%
\usepackage{amsfonts}%
\usepackage{amssymb}%
\usepackage{amsthm}
\usepackage{graphicx}
\usepackage[english]{babel}
\usepackage{csquotes}
\usepackage{tikz-cd}
\newtheorem{Defn}{Definition}
\newtheorem{Prop}{Proposition}

\begin{document}

\begin{frontmatter}

\title{The zero dynamics of feedback linearisation and Ehresmann connections}
\author[rvt]{T.~Chingozha\corref{cor1}}
\author[rvt]{O.T.~Nyandoro}
\author[rvt]{A.~van Wyk}
\author[rvt]{J.E.D.~Ekoru}
\cortext[cor1]{Corresponding author}
\address[rvt]{University of the Witwatersrand, Johannesburg, South Africa}

\begin{abstract}
The notion of zero dynamics is a cornerstone of many solutions to important control problems such as feedback linearisation and disturbance decoupling. For a SISO affine control system with relative degree strictly less than the order of the system, it is known that there exists a state transformation and a feedback transformation such that the system can be transformed to a linear system. Making use of the fibre bundle structure induced by the state transformation we show that it is possible to define a connection on the state manifold such that the zero dynamics can be defined as a vertical vector field. In this formalism the zero dynamics can be understood as motions along the fibres of a fibre bundle.    
\end{abstract}
\begin{keyword}
feedback linearisation \sep Ehresmann connection \sep zero dynamics
\end{keyword}

\end{frontmatter}

\section{Introduction}
The concept of zero dynamics plays a key role in the solution of a number of important control problems for non-linear systems such as asymptotic stabilization and tracking, disturbance decoupling, high gain feedback, feedback linearisation and non-interacting control \cite{Isidori}. In the 1980s the notion of zero dynamics was developed as a non-linear extension of the concept of zeros of a transfer function \cite{Krener}. Consider a SISO affine control system described by the following equation
\begin{eqnarray}
\dot{x} &=& f(x) + g(x)u,\\
y &=& h(x).
\end{eqnarray}
The null observable distribution introduced in \cite{Krener} is defined as the controlled invariant distribution contained in $\text{ker}(dh)$, this distribution enables the partitioning of the non-linear system into observable/unobservable components. From this the zero dynamics are defined as the internal dynamics of the unobservable component of the system. Alternatively \cite{Byrnes} define the zero dynamics as the internal dynamics of the system resulting from the constraint $y(t) = 0$ which is achieved by appropriate choice of control input $u(t)$ and initial conditions. These two definitions of the zero dynamics are not in general equivalent \cite{Alberto}.
The approach taken in this paper follows in the spirit of the definition of \cite{Krener} and it will be shown that the approach presented here leads to a familiar definition of \cite{Krener}. The rest of this paper is organised as follows. Section 2 is a brief presentation of the prerequisite mathematical concepts, Section 3 contains the main result of this paper.

\section{Mathematical Preliminaries}
This section presents a brief survey of the differential geometric concepts that will be used in the main sequel. For a more comprehensive coverage of the concepts presented here the reader may consult \cite{Nomizu}.
The central object in the ensuing analysis is that of a fibre bundle which is defined as follows \cite{SteenRod}:

\begin{Defn}
A fibre bundle is the 4-tuple $(E, M, \pi, F)$ where $E$ and $M$ are smooth manifolds called the \textbf{total space} and \textbf{base space} respectively, $\pi$ is a surjective map from $E$ to $M$, $\pi : E \mapsto M$ which is called the \textbf{projection map} . For some $x$ in $M$ the set $\pi^{-1}(x)$ is called the \textbf{fibre over the point $x$ in M} where $\pi^{-1}(x)$ is required to be homeomorphic to the smooth manifold $F$. $F$ is called the \textbf{typical fibre} over $M$.
\end{Defn}

Intuitively a fibre bundle can be interpreted as a manifold that looks locally like a product of the base space. Consider a smooth vector field on the total space $X \in \Gamma^{\infty}(TE)$ where $\Gamma^{\infty}(TE)$ is the set of all smooth sections of the tangent bundle $TE$ over $E$. The vector field $X$ is called projectable if there exists a vector field on the base space $Y \in \Gamma^{\infty}(TM)$ such that the following holds

\begin{equation}
T_x\pi \circ X(x) = Y\circ \pi(x).
\end{equation}
 
It can easily be shown that if the curve $\gamma(t)$ in $E$ is an an integral curve of the vector field $X$ then the curve $\pi\circ\gamma(t)$ in $M$ is an integral curve of the vector field $Y$. Since $\pi$ is surjective the projection operation is a many to one operation, thus there is no canonical way of relating a vector field (smooth curve) on the base space back to a vector field (smooth curve) on the total space. To be able to uniquely relate a vector field (smooth curve) on the base space back to a vector field (smooth curve) on the total space requires the specification of extra structure and this extra structure comes in the form of a connection.

\begin{Defn}
Let $(M, E, \pi, F)$ be a fibre bundle, the vertical bundle VE is the sub-bundle of the tangent bundle $TE$ defined as
\begin{equation}
VE = \{\mathbf{v}_x \in T_x E | T_x\pi(\mathbf{v}_x) = \mathbf{0}, \forall x \in E\} = \text{ker}(T\pi).
\end{equation}
\end{Defn}

The vertical bundle contains all those tangent vectors that get projected to the zero vector, alternatively the vertical bundle contains all those vectors which are tangent to any fibre. A connection on fibre bundle $(E, M, \pi, F)$ is a sub-bundle of $TE$ which is complementary to the vertical bundle, this sub-bundle is called the horizontal sub-bundle \cite{Cushman}.  This is expressed formally as:

\begin{Defn}
A connection on the fibre bundle $(E,M,\pi, F)$ is a smooth sub-bundle $HE\subset TE$ such that $TE = HE\oplus VE$.
\end{Defn}

The connection makes it possible to lift objects from the base space back to the total space. Consider a vector $Y_p \in T_pM$ on the base space, the horizontal lift $\text{Hor}_x(Y_p) \in T_xE$ is a vector on the total space which gets projected to the vector $Y_p$.  Thus there is a horizontal lift map
\begin{equation}
\text{Hor}_x : T_pM \mapsto T_xE, where \quad \pi(x) = p
\end{equation}

This horizontal lift map satisfies the property that $T_x\pi \circ \text{Hor}_x = \text{id}_{T_pM}$. For a curve $\gamma(t) \in M$ the horizontal lift is curve $\tilde{\gamma}(t) \in E$ such that,

\begin{equation}
\frac{d}{dt}\tilde{\gamma}(t) = \text{Hor}_{\tilde{\gamma}(t)}(\dot{\gamma}(t))
\end{equation}

Consider two vectors on $X, \bar{X} \in T_xE$ which get projected to the same vector i.e $T_x\pi(X) = T_x\pi(\bar{X})$, the difference between these two vectors is a vertical vector ($X - \bar{X} \in V_xE$). 


\section{Feedback Linearisation and Zero Dynamics}
For the purpose of this presentation attention will be restricted to SISO affine control systems. Also the fibre bundle description of Brockett \cite{Brockett} will be used. A control system is defined as the 5-tuple $(E, M, \pi, U, F)$ where $(E, M, \pi, U)$ defines a fibre bundle with typical fibre $U$ and a smooth map $F : E \mapsto TM$ which satisfies $\pi_{TM} \circ F = \pi$ where $\pi_{TM}$ is the tangent bundle projection. Let $\Sigma$ be the affine control system
\begin{eqnarray}
\dot{x} &=& f(x) + g(x)u, \quad x \in M \subset \mathbb{R}^n, u \in U\subset \mathbb{R},\\
y &=& h(x), \quad y\in \mathbb{R}.
\end{eqnarray}

In the fibre bundle representation we have $E = M\times U$ and $F:(x,u) \mapsto f(x) + g(x)u$.  If $\Sigma$ has relative degree $r < n$ there exists a surjective submersion map $\Phi : M \mapsto N\subset \mathbb{R}^r$ defined as
\begin{equation}
\Phi = \left[ h(x), L_fh(x), \cdots,L_f^{r-1}h(x) \right]^T
\end{equation}

and a static feedback $\Psi(x,u) = v$ such that there exists a linear and controllable quotient control system on $N$ such that the following commutes as follows.

\begin{center}
\begin{tikzcd}
E\times\mathbb{R} \arrow{rrr}{(\Phi,\Psi)} \arrow{dd}{\pi_{M\times\mathbb{R}}}\ar{dr}{f}
& & & M\times\mathbb{R} \arrow{dd}{\pi_{M\times\mathbb{R}}} \ar{dl}{\tilde{f}} \\
 & TE \arrow{r}{T\Phi} \ar{dl}{\pi_{TE}} & TM \ar{dr}{\pi_{TM}} & \\
E \arrow{rrr}{\Phi} & & & M
\end{tikzcd}
\end{center}

Since $\Phi$ is a submersion the triple $(E, \Phi, M)$ defines a fibre bundle with total space $E$, base space $M$ and projection map $\Phi$. The fibres of this fibre bundle are the equivalence sets of the submersion map $\Phi$. This insight means the triple can be equipped with a connection $HE$ such that $TE = HE\oplus \text{ker}(T\Phi)$. The connection can be used to define the horizontal lift of the linear controllable system on $S$ as follows:

\begin{Defn}
Given the linear system defined by the 5-tuple $(M\times\mathbb{R}, M, \pi, \mathbb{R}, \tilde{f})$ the horizontal lift of this system is the 5-tuple $(E\times\mathbb{R}, E, \pi, \mathbb{R}, \tilde{f}^H)$ where the map $\tilde{f}^H : E\times\mathbb{R} \mapsto TE$ is defined as

\begin{equation}
\tilde{f}^H(x,u) = \text{Hor}_x\circ\tilde{f}\circ(\Phi, \Psi)(x,u).
\end{equation}
\end{Defn}

To show that the definition of the control system horizontal lift is proper consider the curve $(\gamma^H(t), u^H(t)) \\ \in E\times\mathbb{R}$ such that,

\begin{equation}
\frac{d}{dt}\gamma^H(t) = \text{Hor}_{\gamma^H(t)}\circ \tilde{f} \circ (\Phi, \Psi)\circ (\gamma^H(t), u^H(t)).
\end{equation}
 
That is $\gamma^{H}(t)$ is a trajectory of the horizontally lifted system. Consider the time derivative of the curve $\Phi \circ \gamma^H(t) \in M$,

\begin{eqnarray}
\frac{d}{dt}(\Phi\circ \gamma^H(t)) &=& T_{\gamma^H(t)}\Phi \circ \frac{d}{dt}\gamma^H(t)\\
 &=& T_{\gamma^H(t)}\Phi \circ \text{Hor}_{\gamma^H(t)}\circ \tilde{f} \circ (\Phi, \Psi)\circ (\gamma^H(t), u^H(t))\\
&=& \tilde{f}(\Phi\circ \gamma^H(t), \Psi\circ u(t)).
\end{eqnarray}

This shows that the trajectories of the horizontally lifted system get projected to the original system.

Using the machinery that has been developed so far the zero dynamics are defined as the difference between the original system dynamics and the horizontally lifted linear system dynamics, 

\begin{equation}
f^Z(x,u) = f(x,u) - \text{Hor}_x\circ \tilde{f}\circ (\Phi,\Psi) \circ (x,u).
\end{equation}

By construction the original system dynamics and the horizontally lifted dynamics are horizontal thus the difference between the two is vertical. This means that the zero dynamics are tangent to the fibres (i.e. equivalence sets under $\Phi$ equivalence) of $E$. From this the zero dynamics can be understood as the system motion along the fibres while the linearised dynamics encode the system motion from one fibre to another. If $x_0$ is chosen such that $\Phi(x_0) = 0$ (i.e output is zero)the zero dynamics $f^Z(x_0, u_0)$ coincide with the definition in \cite{Krener} where the zero dynamics are given by the ``internal" dynamics of the system when the input is chosen such that the output stays zero.

\begin{Prop}
Suppose an n-dimensional system has relative degree $r < n$ at $x^0$. Set $\lambda_i = L_f^{i-1}h(x)$ for $i = 1, \cdots, r$, it is always possible to find $n-r$ functions $\lambda_{r+1}(x),\cdots,\eta_n(x)$ such that the mapping
 
\begin{equation*}
\Lambda(x) = \left(
\begin{array}{c}
\lambda_1(x)\\
\cdots\\
\lambda_n(x)
\end{array}
\right)
\end{equation*}

is a diffeomorphism.
\end{Prop}

From the diffeomorphism $\Lambda$ the vertical and horizontal sub-bundles can be defined as $VE = \text{span}([d\lambda_1(x), \\ \cdots,d\lambda_r(x)]^\bot)$ and $HE = \text{span}([d\lambda_{r+1}(x),\cdots, d\lambda_n(x)]^\bot)$. Thus the choice of $n-r$ maps in the above proposition can be understood as a prescription of a connection on the state manifold. To show how these ideas can be applied consider the following example taken from \cite{Isidori}

\subsection{Example}
Consider the SISO system
\begin{eqnarray}
\dot{x} &=& \left(
\begin{array}{c}
exp(x^2)\\
x_1x_2\\
x_2
\end{array} \right)
+ \left(
\begin{array}{c}
exp(x^2)\\
1\\
0
\end{array} \right)u\\
y &=& h(x) = x_3
\end{eqnarray}

For this system the diffeomorphism $\Lambda: \mathbb{R}^3 \mapsto \mathbb{R}^3$ is of the form:
\begin{equation}
\Phi(x) = \left[
\begin{array}{c}
z_1\\
z_2\\
z_3
\end{array} \right] 
=
\left[
\begin{array}{c}
h(x)\\
L_fh(x)\\
\phi_3(x)
\end{array} \right]
=
\left[
\begin{array}{c}
x_3\\
x_2\\
1 + x_1 - exp(x_2)
\end{array} \right]
\end{equation}
It can be easily verified that the map $\Lambda$ is non-singular. From this the submersion map $\Phi$ is given by $\Phi(x) = \left[h(x), L_fh(x) \right]^T$ and the static feedback $\Psi(x, u) = u + x_1x_2$. The horizontal and vertical subspaces are
\begin{eqnarray}
V_x\mathbb{R}^3 &=& \text{span}\{dx_2, dx_3\}^\bot\\
H_x\mathbb{R}^3 &=& \text{span}\{dx_1 - exp(x_2)dx_2\}^\bot
\end{eqnarray}

Projecting the non-linear system to $\mathbb{R}^2$ and applying the feedback transformation produces the familiar feedback linearised dynamics.
\begin{equation}
\left[
\begin{array}{c}
\dot{z}_1\\
\dot{z}_2
\end{array} \right]
=
\left[
\begin{array}{cc}
0 & 1\\
0 & 0
\end{array} \right]
\left[
\begin{array}{c}
z_1\\
z_2
\end{array} \right]
+
\left[
\begin{array}{c}
0\\
1
\end{array} \right] v
\end{equation}

The horizontal lift map is:
\begin{equation}
\text{Hor}_x : T_{(z_1,z_2)}\mathbb{R}^2 \mapsto T_x\mathbb{R}^3
\end{equation}

\begin{equation}
\text{Hor}_x : \left[
\begin{array}{c}
Y^1(z, v)\\
Y^2(z,v)
\end{array} \right]
\mapsto
\left[
\begin{array}{c}
exp(x_2)Y^2(z,v)\\
Y^2(z,v)\\
Y^1(z,v)
\end{array} \right]_{z = \Phi(x), v = u + x_1x_2}
\end{equation}

The lifted dynamics are:
\begin{equation}
Hor_x(\left[
\begin{array}{c}
z_2\\
v
\end{array}\right]) = \left[
\begin{array}{c}
exp(x_2)v\\
v\\
z_2
\end{array} \right]_{z = \Phi(x), v = u + x_1x_2}
=
\left[
\begin{array}{c}
x_1x_2exp(x_2) + uexp(x_2)\\
x_1x_2 + u\\
x_2
\end{array} \right]
\end{equation}
From which the zero dynamics are calculated as:
\begin{equation}
f^Z = \left[
\begin{array}{c}
-x_1(1 + x_2exp(x_2))\\
0\\
0
\end{array}\right]
\end{equation}

The zero dynamics are parameterised by $(x_2, x_3)$ which means that the dynamics vary depending on the fibre however for most control tasks it is required to regulate the output thus by setting $(x_2=0, x_3=0)$ the dynamics become

\begin{equation}
\dot{x}_1 = -x_1.
\end{equation}
 Which agrees with the result in \cite{Isidori}.


\section{Conclusion}
An alternative interpretation of the notion of zero dynamics was presented in this paper. It is shown that, by utilising the fibre bundle structure of the state space manifold, a connection can be defined on the state space manifold. This connection enables the lifting of the linearised dynamics from the base space to the total space. The zero dynamics are then defined as the difference between the original dynamics and the lifted linear dynamics. It then follows from the construction that the zero dynamics represent the motion of the system along the fibres while the linearised dynamics represent the motion of the system from fibre to fibre.

\section*{References}

\bibliographystyle{elsarticle-num}
\bibliography{articleDraft}

\end{document}